# On the Consistency of the Arithmetic System

T. J. Stępień, Ł. T. Stępień

*The Pedagogical University of Cracow, ul. Podchorazych 2, 30 - 084 Krakow, Poland*

**Abstract:** In this paper we establish that the well-known Arithmetic System is consistent in the traditional sense. The proof is done within this Arithmetic System.

## 1. Introduction and Overview

The problem of consistency of the Arithmetic System, is the well-known problem (see [1], [3], [4], [5], [7], [8], [9],[12], [14], [15], [16], [17], [19], [20], [22], [23], [24], [25], [26], [27], [28], [29], [30], [31], [32], [33], [34], [35], [37], [38], [39], [41], [42], [43], [46], [47], [48], [49], [50], [51], [52], [55], [56], [57], [58], [59], [60], [61], [63], [64], [65], [66], [67], [68], [69], [70], [71], [73], [74], [75], [76], [78], [79], [80]). In this paper we give the elementary proof of the consistency of the Arithmetic System. This result was included in the talk, presented under the similar title, at "Logic Colloquium 2009" in Sofia (Bulgaria), [59]. The abstract of this talk has been published in "The Bulletin of Symbolic Logic", [60].

## 2. Terminology

Let: $\rightarrow$, $\sim$, $\vee$, $\wedge$, $\equiv$ denote the connectives of implication, negation, disjunction, conjunction and equivalence, respectively. $\mathcal{N} = \{1,2,\ldots\}$ denotes the set of all natural numbers.

Next $At_0 = \{p_1^1, p_2^1, \ldots, p_1^2, p_2^2, \ldots, p_1^k, p_2^k, \ldots\}$ $(k \in \mathcal{N})$ denotes the set of all propositional variables. Hence, $S_0$ is the set of all well-formed formulas, which are built in the usual manner from propositional variables and by means of logical connectives. $P_0(\phi)$ denotes the set of all propositional variables occuring in $\phi$ $(\phi \in S_0)$. $R_{S_0}$ denotes the set of all rules over $S_0$ (see [44]). $E(\mathfrak{M})$ is the set of all formulas valid in the matrix $\mathfrak{M}$. The $\mathfrak{M}_2$ denotes the classical two-valued matrix. $Z_2$ is the set of all formulas valid in the matrix $\mathfrak{M}_2$ (see [44]).

The symbols $x_1, x_2, \ldots$ are individual variables. $a_1, a_2, \ldots$ are individual constants. $V$ is the set of all individual variables. $C$ is the set of all individual constants. $P_i^n (i, n \in \mathcal{N} = \{1, 2, \ldots\})$ are $n$-ary predicate letters. The symbols $f_i^n (i, n \in \mathcal{N})$ are $n$-ary function letters. The symbols $\bigwedge x_k$, $\bigvee x_k$ are quantifiers. $\bigwedge x_k$ is the universal quantifier and $\bigvee x_k$ is the existential quantifier. The function letters, applied to the individual variables and individual constants, generate terms. The symbols $t_1, t_2, \ldots$ are terms. $T$ is the set of all terms. $V \cup C \subseteq T$. The predicate letters, applied to terms, yield simple formulas, i.e. if $P_i^k$ is a predicate letter and $t_1, \ldots, t_k$ are terms, then $P_i^k(t_1, \ldots, t_k)$ is a simple formula. $Smp$ is the set of all simple formulas. Next, $At_1$ is the set of all atomic formulas, i.e. $At_1 = \{P_i^k(x_{j_1}, \ldots, x_{j_k}): i, k, j_1, \ldots, j_k \in \mathcal{N}\}$. At last, $S_1$ is the set of all well-formed formulas. $FV(\phi)$ denotes the set of all free variables occuring in $\phi$, where $\phi \in S_1$.

$x_k \in Ff(t_m, \phi)$ expresses that $x_k$ is free for term $t_m$ in $\phi$ $(\phi \in S_1)$. By $x_k/t_m$ we denote the substitution of the term $t_m$ for the individual variable $x_k$.

**Corresponding author:** Ł. T. Stępień, The Pedagogical University of Cracow, Kraków, Poland. E-mail: sfstepie@cyf-kr.edu.pl, http://www.ltstepien.up.krakow.pl

$P_1(\phi)$ denotes the set of all predicate letters occuring in $\phi$ ($\phi \in S_1$). If $FV(\phi) = \{x_1, \ldots, x_k\}$, then $\wedge\phi = \wedge x_1 \ldots \wedge x_k \phi$. $R_{S_1}$ denotes the set of all rules over $S_1$. $\overline{S}_1 = \{\phi \in S_1 : FV(\phi) = \emptyset\}$.

We use $\Rightarrow, \neg, \mathbb{V}, \&, \Leftrightarrow, \forall, \exists$ as metalogical symbols. Next, $r_0^i (i \in \{0,1\})$ denotes Modus Ponens for propositional and predicate calculi, respectively.

$r_+$ denotes the generalization rule. $R_{0+} = \{r_0^1, r_+\}$.

$L_2$ is the set of all formulas valid in the classical calculus of quantifiers (see [45]). We write $X \subset Y$ for $X \subseteq Y$ and $Y \neq X$. For any $X \subseteq S_i$, $Cn(R, X)$ is the smallest subset of $S_i$, containing $X$ and closed under the rules $R \subseteq R_{S_i}$ and $i \in \{0,1\}$. The couple $\langle R, X\rangle$ is called a system, whenever $R \subseteq R_{S_i}$ , $X \subseteq S_i$ and $i \in \{0,1\}$ (see [44], [45], *cf.* [11]).

Now we repeat some well-known properties of operation of consequence and some well-known definitions (see [34], [44], [45]). Let $R \subseteq R_{S_i}$ and $X \subseteq S_i$. Then:

$(a_1)$ $X \subseteq Cn(R, X)$,

$(a_2)$ $X \subseteq Y \Rightarrow Cn(R, X) \subseteq Cn(R, Y)$,

$(a_3)$ $R \subseteq R' \Rightarrow Cn(R, X) \subseteq Cn(R', X)$,

$(a_4)$ $Cn(R, Cn(R, X)) \subseteq Cn(R, X)$,

$(a_5)$ $Cn(R, X) = \bigcup\{Cn(R, Y) : Y \in Fin(X)\}$,

where $Y \in Fin(X)$ denotes that $Y$ is the finite subset of $X$ and $i \in \{0,1\}$.

**Definition 1.1.** $\langle R, X\rangle \in Cns^T \Leftrightarrow (\neg\exists\alpha \in S_i)[\alpha \in Cn(R, X)\ \&\ {\sim}\alpha \in Cn(R, X)]$, where $i \in \{0,1\}$.

**Definition 1.2.** $\langle R, X\rangle \in Cns^A \Leftrightarrow Cn(R, X) \neq S_i$, where $i \in \{0,1\}$.

$\langle R, X\rangle \in Cns^T$ denotes that the system $\langle R, X\rangle$ is consistent in the traditional sense, and $\langle R, X\rangle \in Cns^A$ denotes that the system $\langle R, X\rangle$ is consistent in the absolute sense (see [44], [45]).

## 2. Basic Theorems

**Theorem 2.1.** $\langle R_{0+}, L_2\rangle \in Cns^T$.

**Theorem 2.2.** $\langle R_{0+}, L_2\rangle \in Cns^A$.

**Theorem 2.3.** $(\forall\alpha \in \overline{S}_1)(\forall\beta \in S_1)(\forall X \subseteq S_1)$
$[\beta \in Cn(R_{0+}, A_2 \cup X \cup \{\alpha\}) \Rightarrow (\alpha \rightarrow \beta) \in Cn(R_{o+}, A_2 \cup X)]$.

**Theorem 2.4.** $(\forall\alpha \in \overline{S}_1)(\forall X \subseteq S_1)$
$[Cn(R_{0+}, A_2 \cup X \cup \{\alpha\}) = S_1 \Leftrightarrow {\sim}\alpha \in Cn(R_{0+}, A_2 \cup X)]$.

**Theorem 2.5.** $(\forall\alpha \in \overline{S}_1)(\forall X \subseteq S_1)$
$[\alpha \notin Cn(R_{0+}, A_2 \cup X) \Leftrightarrow Cn(R_{0+}, A_2 \cup X \cup \{{\sim}\alpha\}) \neq S_1]$,

where $A_2$ denotes the set of the axioms of the classical calculus of quantifiers (see [45]).

## 3. Arithmetic Terminology

Next, $S_A$ denotes the set of all well-formed formulas of the Arithmetic System. Hence, $Fv(\phi)$ denotes the set of all free variables occuring in $\phi$, where $\phi \in S_A$. $x_k \in Ff_A(t_m, \phi)$ expresses that $x_k$ is free for term $t_m$ in $\phi$, where $\phi \in S_A$.

$\overline{S}_A = \{\phi \in S_A : Fv(\phi) = \emptyset\}$ . Analogically, $R_{S_A}$ denotes the set of all rules over $S_A$. For any $X \subseteq S_A$ and for any $R \subseteq R_{S_A}$, $Cn(R, X)$ is the smallest subset of $S_A$ , containing $X$ and closed under the rules of $R$. The couple $\langle R, X\rangle$ is called a system, whenever $R \subseteq R_{S_A}$ and $X \subseteq S_A$. $R_{0+}^P = \{r_0^P, r_+^P\} \subseteq R_{S_A}$, where $r_0^P$ and $r_+^P$ are Modus Ponens and generalization rule in the Arithmetic System, respectively. By $\psi^1, \psi^2, \psi^3, \psi^4, \psi^5, \psi^6, \psi^7, \psi^8, \psi^9, \psi^{10}, \psi^{11}, \psi^{12}$ , we denote the specific axioms of the Arithmetic System, where:

$\psi^1$. $\wedge x_1 (x_1 = x_1)$,

$\psi^2$. $\wedge x_1 \wedge x_2 (x_1 = x_2 \rightarrow x_2 = x_1)$,

$\psi^3$. $\wedge x_1 \wedge x_2 \wedge x_3 (x_1 = x_2 \rightarrow (x_2 = x_3 \rightarrow x_1 = x_3))$,

$\psi^4$. $\wedge x_1 \wedge x_2 \wedge x_3 \wedge x_4 (x_1 = x_2 \rightarrow (x_3 = x_4 \rightarrow (x_1 + x_3 = x_2 + x_4)))$,

$\psi^5$. $\wedge x_1 \wedge x_2 \wedge x_3 \wedge x_4 (x_1 = x_2 \rightarrow (x_3 = x_4 \rightarrow (x_1 \cdot x_3 = x_2 \cdot x_4)))$,

$\psi^6$. $\wedge x_1 \wedge x_2 \wedge x_3 \wedge x_4 (x_1 = x_2 \rightarrow (x_3 = x_4 \rightarrow (x_1 < x_3 \rightarrow x_2 < x_4)))$,

$\psi^7$. $\wedge x_1\ {\sim}(1 = x_1 + 1)$,

$\psi^8$. $\wedge x_1 \wedge x_2 (x_1 + 1 = x_2 + 1 \rightarrow x_1 = x_2)$,

$\psi^9$. $\wedge x_1 \wedge x_2 (x_1 + (x_2 + 1) = (x_1 + x_2) + 1)$,

$\psi^{10}$. $\wedge x_1 (x_1 \cdot 1 = x_1)$,

$\psi^{11}$. $\wedge x_1 \wedge x_2 [x_1 \cdot (x_2 + 1) = (x_1 \cdot x_2) + x_1]$,

$\psi^{12}$. $\wedge x_1 \wedge x_2 [x_1 < x_2 \equiv \vee x_3 (x_1 + x_3 = x_2)]$.

Hence, $X_P = \{\psi^1, \psi^2, \psi^3, \psi^4, \psi^5, \psi^6, \psi^7, \psi^8, \psi^9, \psi^{10}, \psi^{11}, \psi^{12}\}$. Next, the induction schema is the set of the following axioms:

$$\psi^{13}.\ \Big(\phi(1) \wedge \bigwedge x_1\big(\phi(x_1) \rightarrow \phi(x_1+1)\big)\Big) \rightarrow \bigwedge x_1 \phi(x_1),$$

where $\phi(1), \phi(x), \phi(x+1) \in S_A$.

Hence, $Y_P$ denotes here the set of all axioms of induction. Thus, $L_2^r$ and $A_r$ denote the set of all logical axioms and the set of all specific axioms of the Arithmetic System, respectively, where $A_r = X_P \cup Y_P$ (see [16], [50]).

Hence, $\langle R_{0+}^P, L_2^r \cup A_r\rangle$ is the Arithmetic System. In [50], one can read that the system $\langle R_{0+}^P, L_2^r \cup A_r\rangle$ is a modification of Peano's Arithmetic System.

$Sx$ denotes here the successor of $x$ (see [34]).

Next, by $Q^1, Q^2, Q^3, Q^4, Q^5, Q^6, Q^7, Q^8, Q^9$, we denote the other specific axioms of the arithmetic system,where:

$Q^1.\ \bigwedge x\ (x + 0 = x)$,

$Q^2.\ \bigwedge x \bigwedge y\ (x \cdot Sy = x \cdot y + x)$,

$Q^3.\ \bigwedge x \bigwedge y\ (Sx = Sy \rightarrow x = y)$,

$Q^4.\ \bigwedge x \bigvee y\ (y = Sx)$,

$Q^5.\ \bigwedge x \bigwedge y\ [x + Sy = S(x + y)]$,

$Q^6.\ \bigwedge x\ (x \cdot 0 = 0)$,

$Q^7.\ \sim\bigvee x\ (Sx + 1 = 1)$,

$Q^8.\ \bigwedge x_1 \bigwedge x_2\ [\bigvee x_3\ (Sx_3 + x_1 = x_2) \equiv (x_1 < x_2)]$,

$Q^9.\ \bigwedge x \sim(Sx = 0)$.

Hence,

$X_P' = \{Q^1, Q^2, Q^3, Q^4, Q^5, Q^6, Q^7, Q^8, Q^9\}$.

Next, the induction schema is the set of the following axioms:

$$Q^{10}.\ \Big(\phi(0) \wedge \bigwedge x\big(\phi(x) \rightarrow \phi(Sx)\big)\Big) \rightarrow \bigwedge x \phi(x),$$

where $\phi(1), \phi(x), \phi(Sx) \in S_A$.

Hence, $Y_P'$ denotes the set of all axioms of induction and $A_r'$ denotes the set of all specific axioms of the arithmetic system, where $A_r' = X_P' \cup Y_P'$.

In consequence, one can obtain another arithmetic system (*cf.* [1], [2], [4], [6], [7], [10], [13], [14], [16],[17], [18], [21], [23], [25], [28], [30], [31], [32], [34], [36], [37], [40], [47], [48], [50], [51], [52], [53], [54],[58], [59], [62], [71], [72], [77], [78]), namely, $\langle R_{o+}^P, L_2^{r'} \cup A_r'\rangle$, where $L_2^{r'}$ is the set of all logical axioms.

$L_1^2$ denotes the well-known subset of the set $L_2^r$ (see[16], [44], [45], [50]). Namely, $L_1^2 = \{\phi^1, \phi^2, \phi^3, \phi^4, \phi^5, \phi^6, \phi^7, \phi^8, \phi^9, \phi^{10}, \phi^{11}, \phi^{12}\}$, where:

$\phi^1.\ [\alpha \rightarrow (\beta \rightarrow \gamma)] \rightarrow [(\alpha \rightarrow \beta) \rightarrow (\alpha \rightarrow \gamma)]$,

$\phi^2.\ (\sim\alpha \rightarrow \alpha) \rightarrow \alpha$,

$\phi^3.\ \sim\alpha \rightarrow (\alpha \rightarrow \beta)$,

$\phi^4.\ \alpha \rightarrow (\beta \rightarrow \alpha)$,

$\phi^5.\ \alpha \wedge \beta \rightarrow \alpha$,

$\phi^6.\ \alpha \wedge \beta \rightarrow \beta$,

$\phi^7.\ \alpha \rightarrow (\beta \rightarrow \alpha \wedge \beta)$,

$\phi^8.\ \alpha \rightarrow \alpha \vee \beta$,

$\phi^9.\ \beta \rightarrow \alpha \vee \beta$,

$\phi^{10}.\ (\alpha \rightarrow \beta) \rightarrow [(\delta \rightarrow \beta) \rightarrow (\alpha \vee \delta \rightarrow \beta)]$,

$\phi^{11}.\ \bigwedge x_k \phi \rightarrow \phi\left(\frac{x_k}{t_n}\right)$,if $x_k \in Ff_A(t_n, \phi)$,

$\phi^{12}.\ \bigwedge x_k(\phi \rightarrow \psi) \rightarrow (\phi \rightarrow \bigwedge x_k \psi)$, if $x_k \notin Fv(\phi)$

and

$\alpha, \beta, \gamma, \delta, \phi, \psi \in S_A$.

The analogons of Definition 1.1., Definition 1.2., Theorem 2.1., Theorem 2.2., Theorem 2.4. and Theorem 2.5., are the following (where $R \subseteq R_{S_A}$ and $X \subseteq S_A$):

**Definition 3.1.** $\langle R, X\rangle \in Cns_A^T \Leftrightarrow (\neg\exists\alpha \in S_A)$
$[\alpha \in Cn(R, X)\ \&\ \sim\alpha \in Cn(R, X)]$.

**Definition 3.2.** $\langle R, X\rangle \in Cns_A^A \Leftrightarrow Cn(R, X) \neq S_A$.

**Theorem 3.3.** $\langle R_{0+}^P, L_2^r\rangle \in Cns_A^T$.

**Theorem 3.4.** $\langle R_{0+}^P, L_2^r\rangle \in Cns_A^A$.

**Theorem 3.5.** $(\forall\alpha \in \overline{S}_A)(\forall X \subseteq S_A)$
$[Cn(R_{0+}^P, L_1^2 \cup X \cup \{\alpha\}) = S_A \Leftrightarrow$
$\sim\alpha \in Cn(R_{0+}^P, L_1^2 \cup X)]$.

**Theorem 3.6.** $(\forall\alpha \in \overline{S}_A)(\forall X \subseteq S_A)$
$[\alpha \notin Cn(R_{0+}^P, L_1^2 \cup X) \Leftrightarrow$
$Cn(R_{0+}^P, L_1^2 \cup X \cup \{\sim\alpha\}) \neq S_A]$,

where $L_1^2 \subseteq L_2^r$ (see [11], [16], [34], [45], [50], [65]).

## 4. The Basic Corollaries and Lemmas

At first, we introduce the following formulas:

$(I_1)\ O_0 = \psi^7 \equiv \sim\sim\psi^1$,

$(I_2)\ u_{27} = \sim(1 < 1)$,

$(I_3)\ O_6 = O_0 \to (\psi^7 \to \psi^1)$,

$(I_4)\ \alpha_2^x = \psi^1 \to \psi^7$,

$(I_5)\ \gamma_2^{'} = O_0 \to u_{27}$,

$(I_6)\ \gamma_0^{'} = (\psi^7 \to \psi^1) \to \psi^{12}$,

$(I_7)\ \gamma_0 = u_{27} \to \gamma_0^{'}$,

$(I_8)\ \gamma_4^{'} = \gamma_0^{'} \to O_0$.

Next, we assume that

$(I_9)\ (\forall\alpha,\delta \in S_A)[\,_\alpha\delta = \alpha \to \delta]$,

$(I_{10})\ (\forall X \subseteq S_A)(\forall\alpha \in S_A)[\,_\alpha X = (\alpha \to \beta{:}\, \beta \in X)]$.

Next, we define the sets $L_1^1$ and $_{\psi^7 O_0 u_{27}\sim\psi^1}L_2^r$, as follows:

$(I_{11})\ L_1^1 = L_1^2 \cup \{\psi^1 \to (\psi^7 \to (\psi^{12} \to \omega)){:}\ \omega \in L_2^r - L_1^2\}$,

$(I_{12})\ _{\psi^7 O_0 u_{27}\sim\psi^1}L_2^r = \{\psi^7 \to (O_0 \to (u_{27} \to (\sim\psi^1 \to \beta))){:}\, \beta \in L_2^r - L_1^2\}$.

**Lemma 4.1.** $\langle R_{0+}^P, L_2^r \cup \{\psi^1,\psi^7,\psi^{12}\}\rangle \notin Cns_A^T \Rightarrow (\forall\alpha \in \overline{S}_A - A^0)\left(\forall\overset{00}{\delta} \in \overline{S}_A - A^0\right)$
$[Cn(R_{0+}^P, L_1^1 \cup {}_{\psi^7 O_0 u_{27}\sim\psi^1}L_2^r \cup N\dot{\overline{\Psi}}_3 \cup \left\{\overset{00}{\delta}\right\} \cup \{\alpha \to \psi^7\} \cup \{\xi\}) = S_A \Rightarrow Cn(R_{0+}^P, L_1^1 \cup {}_{\psi^7 O_0 u_{27}\sim\psi^1}L_2^r \cup N\dot{\overline{\Psi}}_3 \cup \{\alpha \to \psi^7\} \cup \{\xi\}) = S_A]$,

where

$\overset{00}{\delta} = \psi^7 \to \overset{0}{\delta}, \xi = \psi^7 \to (\psi^1 \to \psi^{12})$,

$A^0 = Cn\left(R_{0+}^P, L_1^1 \cup {}_{\psi^7 O_0 u_{27}\sim\psi^1}L_2^r \cup N\dot{\overline{\Psi}}_3 \cup \{\xi\}\right)$,

$N\dot{\overline{\Psi}}_3 = \{O_0 \to \gamma_0, \gamma_2^{'}, \gamma_4^{'}\}$.

Proof. Let

(1) $\langle R_{0+}^P, L_2^r \cup \{\psi^1,\psi^7,\psi^{12}\}\rangle \notin Cns_A^T$

and

(2) $\neg(\forall\alpha \in \overline{S}_A - A^0)\left(\forall\overset{00}{\delta} \in \overline{S}_A - A^0\right)$
$[Cn(R_{0+}^P, L_1^1 \cup {}_{\psi^7 O_0 u_{27}\sim\psi^1}L_2^r \cup N\dot{\overline{\Psi}}_3 \cup \left\{\overset{00}{\delta}\right\} \cup \{\alpha \to \psi^7, \xi\}) = S_A \Rightarrow$
$Cn(R_{0+}^P, L_1^1 \cup {}_{\psi^7 O_0 u_{27}\sim\psi^1}L_2^r \cup N\dot{\overline{\Psi}}_3 \cup \{\alpha \to \psi^7, \xi\}) = S_A]$,

where

(3) $\overset{00}{\delta} = \psi^7 \to \overset{0}{\delta}$,

(4) $\xi = \psi^7 \to (\psi^1 \to \psi^{12})$,

(5) $A^0 = Cn(R_{0+}^P, L_1^1 \cup {}_{\psi^7 O_0 u_{27}\sim\psi^1}L_2^r \cup N\dot{\overline{\Psi}}_3 \cup \{\xi\})$,

(6) $N\dot{\overline{\Psi}}_3 = \{O_0 \to \gamma_0, \gamma_2^{'}, \gamma_4^{'}\}$.

Hence, it follows that

(7) $(\exists\alpha^{'} \in \overline{S}_A - A^0)\left(\exists\overset{00}{\delta}{}^{'} \in \overline{S}_A - A^0\right)$
$[Cn(R_{0+}^P, L_1^1 \cup {}_{\psi^7 O_0 u_{27}\sim\psi^1}L_2^r \cup N\dot{\overline{\Psi}}_3 \cup \left\{\overset{00}{\delta}{}^{'}\right\} \cup \{\alpha^{'} \to \psi^7, \xi\}) = S_A\ \&\ Cn(R_{0+}^P, L_1^1 \cup {}_{\psi^7 O_0 u_{27}\sim\psi^1}L_2^r \cup N\dot{\overline{\Psi}}_3 \cup \{\alpha^{'} \to \psi^7, \xi\}) = A \neq S_A]$,

where

(8) $\overset{00}{\delta}{}^{'} = \psi^7 \to \overset{0}{\delta}{}^{'}$,

(9) $\xi = \psi^7 \to (\psi^1 \to \psi^{12})$,

(10) $A^0 = Cn(R_{0+}^P, L_1^1 \cup {}_{\psi^7 O_0 u_{27}\sim\psi^1}L_2^r \cup N\dot{\overline{\Psi}}_3 \cup \{\xi\})$,

(11) $N\dot{\overline{\Psi}}_3 = \{O_0 \to \gamma_0, \gamma_2^{'}, \gamma_4^{'}\}$.

From (7) – (11), it follows that

(12) $\left(\exists\alpha' \in \overline{S}_A - A^0\right)\left(\exists \overset{00}{\delta}{}' \in \overline{S}_A - A^0\right)$

$[Cn(R^P_{0+}, L^1_1 \cup {}_{\psi^7 O_0 u_{27} \sim\psi^1}L^r_2 \cup N\dot{\overline{\Psi}}_3 \cup$

$\{\alpha' \to \psi^7, \xi\}) = A \;\&\; A \neq S_A \;\&\; \overset{00}{\delta}{}' \notin A].$

From (12), by **Theorem 3.6.**, one can obtain that

(13) $\left(\exists\alpha' \in \overline{S}_A - A^0\right)\left(\exists \overset{00}{\delta}{}' \in \overline{S}_A - A^0\right)$

$[Cn(R^P_{0+}, L^1_1 \cup {}_{\psi^7 O_0 u_{27} \sim\psi^1}L^r_2 \cup N\dot{\overline{\Psi}}_3 \cup$

$\{\alpha' \to \psi^7, \xi\}) = A \;\&\; Cn(R^P_{0+}, L^1_1 \cup$

${}_{\psi^7 O_0 u_{27} \sim\psi^1}L^r_2 \cup N\dot{\overline{\Psi}}_3 \cup \left\{\sim\overset{00}{\delta}{}'\right\} \cup$

$\{\alpha' \to \psi^7\} \cup \{\xi\}) = A^* \;\&\; A^* \neq S_A].$

Hence, from (1), (8), (9), (11), $(I_5)$, $(I_6)$, $(I_8)$, $(I_{11})$, $(I_{12})$ and **Theorem 3.5.**, it follows that

(14) $\left(\exists\alpha' \in \overline{S}_A - A^0\right)\left(\exists \overset{00}{\delta}{}' \in \overline{S}_A - A^0\right)$

$[Cn(R^P_{0+}, L^1_1 \cup {}_{\psi^7 O_0 u_{27} \sim\psi^1}L^r_2 \cup N\dot{\overline{\Psi}}_3 \cup$

$\left\{\sim\overset{00}{\delta}{}'\right\} \cup \{\alpha' \to \psi^7\} \cup \{\xi\}) = A^* \;\&$

$\psi^7, \gamma'_4, \gamma'_2, O_0 \to \gamma_0, \psi^1 \to \psi^{12}, \gamma'_0,$

$O_0, u_{27}, \psi^{12} \to (\psi^7 \to \sim\psi^1), \psi^{12} \to \sim\psi^1,$

$\sim\psi^1 \in A^* \;\&\; A^* \neq S_A].$

Hence and from $(I_{12})$, it follows that

(15) $\left(\exists\alpha' \in \overline{S}_A - A^0\right)\left(\exists \overset{00}{\delta}{}' \in \overline{S}_A - A^0\right)$

$[Cn(R^P_{0+}, L^1_1 \cup {}_{\psi^7 O_0 u_{27} \sim\psi^1}L^r_2 \cup N\dot{\overline{\Psi}}_3 \cup$

$\left\{\sim\overset{00}{\delta}{}'\right\} \cup \{\alpha' \to \psi^7\} \cup \{\xi\}) = A^* \;\&$

$L^r_2 \subseteq A^* \;\&\; A^* \neq S_A].$

From (14), (15) and $(I_1)$, it follows that

(16) $\left(\exists\alpha' \in \overline{S}_A - A^0\right)\left(\exists \overset{00}{\delta}{}' \in \overline{S}_A - A^0\right)$

$[Cn(R^P_{0+}, L^1_1 \cup {}_{\psi^7 O_0 u_{27} \sim\psi^1}L^r_2 \cup N\dot{\overline{\Psi}}_3 \cup$

$\left\{\sim\overset{00}{\delta}{}'\right\} \cup \{\alpha' \to \psi^7\} \cup \{\xi\}) = A^* \;\&$

$\psi^7, \psi^7 \to \psi^1, \psi^1, \sim\psi^1 \in A^* \;\&$

$A^* \neq S_A \;\&\; A^* = S_A].$

Hence, it follows that

(17) $\left(\exists\alpha' \in \overline{S}_A - A^0\right)\left(\exists \overset{00}{\delta}{}' \in \overline{S}_A - A^0\right)$

$[Cn(R^P_{0+}, L^1_1 \cup {}_{\psi^7 O_0 u_{27} \sim\psi^1}L^r_2 \cup N\dot{\overline{\Psi}}_3 \cup$

$\left\{\sim\overset{00}{\delta}{}'\right\} \cup \{\alpha' \to \psi^7\} \cup \{\xi\}) = A^* \;\&$

$A^* \neq S_A \;\&\; A^* = S_A].$

Contradiction □

**Lemma 4.2.** $\langle R^P_{0+}, L^r_2 \cup \{\psi^1, \psi^7, \psi^{12}\}\rangle \notin Cns^T_A \Rightarrow$ $\left(\exists\alpha' \in \overline{S}_A - A^0\right)\left(\forall \overset{00}{\delta} \in \overline{S}_A - A^0\right)[Cn(R^P_{0+}, L^1_1 \cup$ ${}_{\psi^7 O_0 u_{27} \sim\psi^1}L^r_2 \cup N\dot{\overline{\Psi}}_3 \cup \left\{\overset{00}{\delta}\right\} \cup \{\alpha' \to \psi^7, \xi\}) = S_A],$

where

$\overset{00}{\delta} = \psi^7 \to \overset{0}{\delta}, \xi = \psi^7 \to (\psi^1 \to \psi^{12}),$

$A^0 = Cn\left(R^P_{0+}, L^1_1 \cup {}_{\psi^7 O_0 u_{27} \sim\psi^1}L^r_2 \cup N\dot{\overline{\Psi}}_3 \cup \{\xi\}\right),$

$N\dot{\overline{\Psi}}_3 = \{O_0 \to \gamma_0, \gamma'_2, \gamma'_4\}.$

Proof. Suppose to the contrary that

(1) $\langle R^P_{0+}, L^r_2 \cup \{\psi^1, \psi^7, \psi^{12}\}\rangle \notin Cns^T_A$

and

(2) $\neg\left(\exists\alpha' \in \overline{S}_A - A^0\right)\left(\forall \overset{00}{\delta} \in \overline{S}_A - A^0\right)$

$[Cn(R^P_{0+}, L^1_1 \cup {}_{\psi^7 O_0 u_{27} \sim\psi^1}L^r_2 \cup N\dot{\overline{\Psi}}_3 \cup \left\{\overset{00}{\delta}\right\} \cup$

$\{\alpha' \to \psi^7, \xi\}) = S_A],$

where

(3) $\overset{00}{\delta} = \psi^7 \to \overset{0}{\delta},$

(4) $\xi = \psi^7 \to (\psi^1 \to \psi^{12}),$

(5) $A^0 = Cn(R^P_{0+}, L^1_1 \cup {}_{\psi^7 O_0 u_{27} \sim\psi^1}L^r_2 \cup$

$N\dot{\overline{\Psi}}_3 \cup \{\xi\}),$

(6) $N\dot{\overline{\Psi}}_3 = \{O_0 \to \gamma_0, \gamma'_2, \gamma'_4\}.$

From (2) – (6), it follows that

(7) $(\forall\alpha \in \overline{S}_A - A^0)\left(\exists \overset{00}{\delta}' \in \overline{S}_A - A^0\right)$

$[Cn(R^P_{0+}, L^1_1 \cup {}_{\psi^7 O_0 u_{27} \sim \psi^1}L^r_2 \cup N\dot{\overline{\Psi}}_3 \cup \left\{\overset{00}{\delta}'\right\} \cup$

$\{\alpha \to \psi^7, \xi\}) = A^+ \neq S_A]$,

where

(8) $\overset{00}{\delta}' = \psi^7 \to \overset{0}{\delta}'$,

(9) $\xi = \psi^7 \to (\psi^1 \to \psi^{12})$,

(10) $A^0 = Cn(R^P_{0+}, L^1_1 \cup {}_{\psi^7 O_0 u_{27} \sim \psi^1}L^r_2 \cup$

$N\dot{\overline{\Psi}}_3 \cup \{\xi\})$,

(11) $N\dot{\overline{\Psi}}_3 = \{O_0 \to \gamma_0, \gamma_2', \gamma_4'\}$.

From (10) and $(I_{11})$, it follows that

(12) $\sim\psi^7 \in A^0 \Rightarrow \left(\forall \overset{0}{\delta} \in \overline{S}_A\right)\left[\psi^7 \to \overset{0}{\delta} \in A^0\right]$.

Hence, from (8), it follows that

(13) $\sim\psi^7 \in A^0 \Rightarrow \overset{00}{\delta}' \in A^0$.

Hence, from (7) – (11), it follows that

(14) $\sim\psi^7 \notin A^0$.

From (1), (7) – (11), $(I_5), (I_6), (I_{11}), (I_{12})$ and **Theorem 3.5.**, it follows that

(15) $(\forall\alpha \in \overline{S}_A - A^0)\left(\exists \overset{00}{\delta}' \in \overline{S}_A - A^0\right)$

$[Cn(R^P_{0+}, L^1_1 \cup {}_{\psi^7 O_0 u_{27} \sim \psi^1}L^r_2 \cup N\dot{\overline{\Psi}}_3 \cup \left\{\overset{00}{\delta}'\right\} \cup$

$\{\alpha \to \psi^7, \xi\}) = A^+$ & $A^+ \neq S_A$ &

$O_0 \to \gamma_0, \gamma_2', \gamma_4', \psi^1 \to (\psi^7 \to \psi^{12})$,

$\psi^{12} \to (\psi^7 \to \sim\psi^1), \psi^7 \to \sim\psi^1, \gamma_0' \to O_0$,

$\psi^7 \to (\sim\psi^1 \to O_0), \psi^7 \to O_0, O_0 \to \gamma_0'$,

$\psi^7 \to \gamma_0', \psi^7 \to u_{27} \in A^+]$.

From (15), $(I_5)$, it follows that

(16) $(\forall\alpha \in \overline{S}_A - A^0)\left(\exists \overset{00}{\delta}' \in \overline{S}_A - A^0\right)$

$\left[{}_{\psi^7 O_0 \sim \psi^1}L^r_2 \subseteq A^+\right]$.

Hence, from (15), it follows that

(17) $(\forall\alpha \in \overline{S}_A - A^0)\left(\exists \overset{00}{\delta}' \in \overline{S}_A - A^0\right)[O_6 \in A^+]$,

where

(18) $O_6 = O_0 \to (\psi^7 \to \psi^1)$.

Hence, from (15), $(I_{11})$, it follows that

(19) $(\forall\alpha \in \overline{S}_A - A^0)\left(\exists \overset{00}{\delta}' \in \overline{S}_A - A^0\right)$

$[Cn(R^P_{0+}, L^1_1 \cup {}_{\psi^7 O_0 u_{27} \sim \psi^1}L^r_2 \cup N\dot{\overline{\Psi}}_3 \cup$

$\left\{\overset{00}{\delta}'\right\} \cup \{\alpha \to \psi^7, \xi\}) = A^+$ & $A^+ \neq S_A$ &

$\gamma_0' \to (\psi^7 \to \psi^1), \psi^7 \to \psi^1, \psi^1 \to \sim\psi^7$,

$\sim\psi^7 \in A^+]$.

Hence, from (10), it follows that

(20) $(\forall\alpha \in \overline{S}_A - A^0)\left(\exists \overset{00}{\delta}' \in \overline{S}_A - A^0\right)$

$[Cn(R^P_{0+}, L^1_1 \cup {}_{\psi^7 O_0 u_{27} \sim \psi^1}L^r_2 \cup N\dot{\overline{\Psi}}_3 \cup \left\{\overset{00}{\delta}'\right\} \cup$

$\{\alpha \to \psi^7, \xi\}) = A^+$ & $A^+ \neq S_A$ &

$A^+ - A^0 = \emptyset$ & $A^0 \subseteq A^+]$.

Hence, from (10) and (19), it follows that

(21) $\sim\psi^7 \in A^0$,

what contradicts with (14). □

**Theorem 4.1**. $\langle R^P_{0+}, L^r_2 \cup \{\psi^1, \psi^7, \psi^{12}\}\rangle \in Cns^T_A$.

Proof. Let

(1) $\langle R^P_{0+}, L^r_2 \cup \{\psi^1, \psi^7, \psi^{12}\}\rangle \notin Cns^T_A$.

From (1), by **Lemma 4.1.** and **Lemma 4.2.**, it follows that

(2) $(\exists\alpha' \in \overline{S}_A - A^0)[Cn(R^P_{0+}, L^1_1 \cup {}_{\psi^7 O_0 u_{27} \sim \psi^1}L^r_2 \cup$

$N\dot{\overline{\Psi}}_3 \cup \{\alpha' \to \psi^7\} \cup \{\xi\}) = S_A]$,

where

(3) $A^0 = Cn(R^P_{0+}, L^1_1 \cup {}_{\psi^7 O_0 u_{27} \sim \psi^1}L^r_2 \cup$

$N\dot{\overline{\Psi}}_3 \cup \{\xi\})$.

Hence, by **Theorem 3.5.**, it follows that

(4) $(\exists\alpha' \in \overline{S}_A - A^0)[\sim(\alpha' \to \psi^7) \in A^0]$,

where

(5) $A^0 = Cn(R^P_{0+}, L^1_1 \cup {}_{\psi^7 O_0 u_{27} \sim \psi^1} L^r_2 \cup N\dot{\bar{\Psi}}_3 \cup \{\xi\})$.

Hence, it follows that

(6) $(\exists \alpha' \in \overline{S}_A)[\alpha' \notin A^0 \ \& \ \alpha' \in A^0]$,

where

(7) $A^0 = Cn(R^P_{0+}, L^1_1 \cup {}_{\psi^7 O_0 u_{27} \sim \psi^1} L^r_2 \cup N\dot{\bar{\Psi}}_3 \cup \{\xi\})$.

Contradiction □

Now,

**Corollary 4.1.** $\langle R^P_{0+}, L^r_2 \cup A_r \rangle \notin Cns^T_A \Rightarrow (\exists Y''_P \subseteq Y_P)(\exists X''_p \subseteq X_P)\{Cn(R^P_{0+}, L^r_2 \cup Y''_p \cup X''_P) = S_A \ \& \ (\forall Y \subset X''_P \cup Y''_P)[Cn(R^P_{0+}, L^r_2 \cup Y) \neq S_A]\}$,

where

$Y''_P = \{\alpha_1, \ldots, \alpha_k\}, X''_p = \{\alpha'_1, \ldots, \alpha'_n\}$ and $k, n \in \mathcal{N}$ and $X_P \cup Y_P = A_r$ .

Proof. By **Theorem 3.3. – Theorem 3.6.**, by the definition of operation of consequence, and by the definitions of the sets $L^r_2$ and $A_r$. □

**Corollary 4.2.** $\langle R^P_{0+}, L^r_2 \cup A_r \rangle \notin Cns^T_A \Rightarrow (\exists Y''_P \subseteq Y_P)$
$(\exists X''_p \subseteq X_P)[Cn(R^P_{0+}, L^r_2 \cup Y''_P \cup X''_P) = S_A \ \&$
$(\forall Y \subset X''_P \cup Y''_P)[Cn(R^P_{0+}, L^r_2 \cup Y) \neq S_A \ \&$
$(\exists Z_P \subseteq \{\psi^2, \psi^3, \psi^4, \psi^5, \psi^6, \psi^8, \psi^9, \psi^{10}, \psi^{11}\} \cup Y_P)$
$[Z_P \subseteq X''_P \cup Y''_P \ \& \ Z_P \neq \emptyset]]]$,

where

$Y''_P = \{\alpha_1, \ldots, \alpha_k\}, X''_P = \{\alpha'_1, \ldots, \alpha'_n\}$ and $k, n \in \mathcal{N}$ and $Y_P \cup X_P = A_r$.

Proof. By **Corollary 4.1.** and by **Theorem 4.1.** □

Using **Corollary 4.2.**, we define the formulas $\overset{0}{\beta}$, $\overset{1}{\beta}$, as follows:

$(I_{14})$ $\overset{0}{\beta} = \alpha_1 \wedge \ldots \wedge \alpha_k \wedge \alpha'_1 \wedge \ldots \wedge \alpha'_n \wedge \psi^1 \wedge \psi^7 \wedge \psi^{12}$,

$(I_{15})$ $\overset{1}{\beta} = \psi^1 \to \left(\psi^7 \to \left(\psi^{12} \to \overset{0}{\beta}\right)\right)$,

where

$\{\alpha_1, \ldots, \alpha_k\} = Y''_P$ and $\{\alpha'_1, \ldots, \alpha'_n\} = X''_P$ and $k, n \in \mathcal{N}$.

Next, we define some sets $L^1_T$ and $N\ddot{\bar{\Psi}}_3$, as follows:

$(I_{16})$ $L^1_T = L^2_1 \cup \left\{\overset{0}{\beta} \to \omega : \omega \in L^r_2 - L^2_1\right\}$,

$(I_{17})$ $N\ddot{\bar{\Psi}}_3 = \left\{O_0 \to \left(u_{27} \to \overset{1}{\beta}\right), O_0 \to \gamma_0, \gamma'_2, \gamma'_4\right\}$.

Thus,

**Corollary 4.3.** $\langle R^P_{0+}, L^r_2 \cup A_r \rangle \notin Cns^T_A \Rightarrow$

$Cn\left(R^P_{0+}, L^1_T \cup \left\{\overset{1}{\beta}\right\} \cup \{\psi^1, \psi^7, \psi^{12}\}\right) = S_A$.

Proof. From $(I_{14})$, $(I_{15})$, $(I_{16})$ and by **Corollary 4.2.** and by the definition of the formula $\overset{1}{\beta}$. □

**Corollary 4.4.** $\langle R^P_{0+}, L^r_2 \cup A_r \rangle \notin Cns^T_A \Rightarrow$

$\left[\sim \overset{0}{\beta} \in Cn(R^P_{0+}, L^1_T)\right]$.

Proof. Let (1) $\langle R^P_{0+}, L^r_2 \cup A_r \rangle \notin Cns^T_A$ and (2) $\sim \overset{0}{\beta} \notin Cn(R^P_{0+}, L^1_T)$. Hence, by **Theorem 3.6.** and $(I_{16})$, we get that (3) $Cn\left(R^P_{0+}, L^1_T \cup \left\{\overset{0}{\beta}\right\}\right) = A \neq S_A$. Hence, from (1), $(I_{14}), (I_{15})$, $(I_{16})$ and by **Corollary 4.3.**, one can get that (4) $A = S_A$, what contradicts (3). □

Now,

**Lemma 4.3.** $\langle R^P_{0+}, L^r_2 \cup A_r \rangle \notin Cns^T_A \Rightarrow$

$(\forall \alpha \in \overline{S}_A - A^1)\left(\forall \overset{00}{\delta} \in \overline{S}_A - A^1\right)$

$[Cn(R^P_{0+}, L^1_1 \cup {}_{\psi^7 O_0 u_{27} \sim \psi^1} L^r_2 \cup N\ddot{\bar{\Psi}}_3 \cup \left\{\overset{00}{\delta}\right\} \cup \{\alpha \to \psi^7, \xi\}) = S_A \Rightarrow Cn(R^P_{0+}, L^1_1 \cup {}_{\psi^7 O_0 u_{27} \sim \psi^1} L^r_2 \cup N\ddot{\bar{\Psi}}_3 \cup \{\alpha \to \psi^7, \xi\}) = S_A]$,

where

$\overset{00}{\delta} = \psi^7 \to \overset{0}{\delta}, \xi = \psi^7 \to (\psi^1 \to \psi^{12})$,

$A^1 = Cn\left(R^P_{0+}, L^1_1 \cup {}_{\psi^7 O_0 u_{27} \sim \psi^1} L^r_2 \cup N\ddot{\bar{\Psi}}_3 \cup \{\xi\}\right)$.

Proof. Let

(1) $\langle R^P_{0+}, L^r_2 \cup A_r \rangle \notin Cns^T_A$

and

(2) $\neg\left(\forall\alpha \in \overline{S}_A - A^1\right)\left(\forall \overset{00}{\delta} \in \overline{S}_A - A^1\right)$

$[Cn(R_{0+}^P, L_1^1 \cup {}_{\psi^7 O_0 u_{27} \sim \psi^1}L_2^r \cup N\ddot{\overline{\Psi}}_3 \cup \left\{\overset{00}{\delta}\right\} \cup \{\alpha \to \psi^7, \xi\}) = S_A \Rightarrow Cn(R_{0+}^P, L_1^1 \cup {}_{\psi^7 O_0 u_{27} \sim \psi^1}L_2^r \cup N\ddot{\overline{\Psi}}_3 \cup \{\alpha \to \psi^7, \xi\}) = S_A]$,

where

(3) $\overset{00}{\delta} = \psi^7 \to \overset{0}{\delta}$,

(4) $\xi = \psi^7 \to (\psi^1 \to \psi^{12})$,

(5) $A^1 = Cn(R_{0+}^P, L_1^1 \cup {}_{\psi^7 O_0 u_{27} \sim \psi^1}L_2^r \cup N\ddot{\overline{\Psi}}_3 \cup \{\xi\})$.

From (2) – (5), it follows that

(6) $\left(\exists\alpha' \in \overline{S}_A - A^1\right)\left(\exists \overset{00}{\delta}{}' \in \overline{S}_A - A^1\right)$

$[Cn(R_{0+}^P, L_1^1 \cup {}_{\psi^7 O_0 u_{27} \sim \psi^1}L_2^r \cup N\ddot{\overline{\Psi}}_3 \cup \left\{\overset{00}{\delta}{}'\right\} \cup \{\alpha' \to \psi^7, \xi\}) = S_A \;\&\; Cn(R_{0+}^P, L_1^1 \cup {}_{\psi^7 O_0 u_{27} \sim \psi^1}L_2^r \cup N\ddot{\overline{\Psi}}_3 \cup \{\alpha' \to \psi^7, \xi\}) = A \neq S_A]$, where

(7) $\overset{00}{\delta}{}' = \psi^7 \to \overset{0}{\delta}{}'$,

(8) $\xi = \psi^7 \to (\psi^1 \to \psi^{12})$,

(9) $A^1 = Cn(R_{0+}^P, L_1^1 \cup {}_{\psi^7 O_0 u_{27} \sim \psi^1}L_2^r \cup N\ddot{\overline{\Psi}}_3 \cup \{\xi\})$.

From (6) – (9), it follows that

(10) $\left(\exists\alpha' \in \overline{S}_A - A^1\right)\left(\exists \overset{00}{\delta}{}' \in \overline{S}_A - A^1\right)$

$[Cn(R_{0+}^P, L_1^1 \cup {}_{\psi^7 O_0 u_{27} \sim \psi^1}L_2^r \cup N\ddot{\overline{\Psi}}_3 \cup \{\alpha' \to \psi^7, \xi\}) = A \;\&\; A \neq S_A \;\&\; \overset{00}{\delta}{}' \notin A]$.

From (10) and by **Theorem 3.6.**, one can obtain that

(11) $\left(\exists\alpha' \in \overline{S}_A - A^1\right)\left(\exists \overset{00}{\delta}{}' \in \overline{S}_A - A^1\right)$

$[Cn(R_{0+}^P, L_1^1 \cup {}_{\psi^7 O_0 u_{27} \sim \psi^1}L_2^r \cup N\ddot{\overline{\Psi}}_3 \cup \{\alpha' \to \psi^7, \xi\}) = A \;\&\; Cn(R_{0+}^P, L_1^1 \cup {}_{\psi^7 O_0 u_{27} \sim \psi^1}L_2^r \cup N\ddot{\overline{\Psi}}_3 \cup \left\{\sim\overset{00}{\delta}{}'\right\} \cup \{\alpha' \to \psi^7\} \cup \{\xi\}) = A^* \;\&\; A^* \neq S_A]$.

From $(I_{11})$, $(I_{14})$, $(I_{16})$, it follows that

(12) $Cn(R_{0+}^P, L_T^1) \subseteq Cn(R_{0+}^P, L_1^1)$.

From (1) – (12), $(I_5)$, $(I_6)$, $(I_8)$, $(I_{14})$ – $(I_{17})$, by **Corollary 4.4.**, it follows that

(13) $\left(\exists\alpha' \in \overline{S}_A - A^1\right)\left(\exists \overset{00}{\delta}{}' \in \overline{S}_A - A^1\right)$

$[Cn(R_{0+}^P, L_1^1 \cup {}_{\psi^7 O_0 u_{27} \sim \psi^1}L_2^r \cup N\ddot{\overline{\Psi}}_3 \cup \left\{\sim\overset{00}{\delta}{}'\right\} \cup \{\alpha' \to \psi^7\} \cup \{\xi\}) = A^* \;\&\; \psi^7, \gamma_4', \gamma_2', O_0 \to \gamma_0, \psi^1 \to \psi^{12}, \gamma_0', O_0, u_{27}, \overset{1}{\beta}, \psi^1 \to \overset{0}{\beta}, \sim\overset{0}{\beta}, \sim\psi^1 \in A^* \;\&\; A^* \neq S_A]$.

Hence, from $(I_{12})$, it follows that

(14) $\left(\exists\alpha' \in \overline{S}_A - A^1\right)\left(\exists \overset{00}{\delta}{}' \in \overline{S}_A - A^1\right)$

$[Cn(R_{0+}^P, L_1^1 \cup {}_{\psi^7 O_0 u_{27} \sim \psi^1}L_2^r \cup N\ddot{\overline{\Psi}}_3 \cup \left\{\sim\overset{00}{\delta}{}'\right\} \cup \{\alpha' \to \psi^7\} \cup \{\xi\}) = A^* \;\&\; L_2^r \subseteq A^* \;\&\; A^* \neq S_A]$.

From (13), (14), and $(I_1)$, it follows that

(15) $\left(\exists\alpha' \in \overline{S}_A - A^1\right)\left(\exists \overset{00}{\delta}{}' \in \overline{S}_A - A^1\right)$

$[Cn(R_{0+}^P, L_1^1 \cup {}_{\psi^7 O_0 u_{27} \sim \psi^1}L_2^r \cup N\ddot{\overline{\Psi}}_3 \cup \left\{\sim\overset{00}{\delta}{}'\right\} \cup \{\alpha' \to \psi^7\} \cup \{\xi\}) = A^* \;\&\; \psi^7, \psi^7 \to \psi^1, \sim\psi^1, \psi^1 \in A^* \;\&\; A^* \neq S_A]$.

From (14) and (15), it follows that

(16) $(\exists\alpha' \in \overline{S}_A - A^1)\left(\exists \overset{00}{\delta}{}' \in \overline{S}_A - A^1\right)$

$[Cn(R^P_{0+}, L^1_1 \cup {}_{\psi^7 O_0 u_{27} \sim\psi^1}L^r_2 \cup N\ddot{\overline{\Psi}}_3 \cup \left\{\sim\overset{00}{\delta}{}'\right\} \cup \{\alpha' \to \psi^7\} \cup \{\xi\}) = A^* \ \&$

$A^* \neq S_A \ \& \ A^* = S_A]$.

Contradiction. □

**Lemma 4.4.** $\langle R^P_{0+}, L^r_2 \cup A_r\rangle \notin Cns^T_A \Rightarrow$

$(\exists\alpha' \in \overline{S}_A - A^1)\left(\forall \overset{00}{\delta} \in \overline{S}_A - A^1\right)[Cn(R^P_{0+}, L^1_1 \cup {}_{\psi^7 O_0 u_{27} \sim\psi^1}L^r_2 \cup N\ddot{\overline{\Psi}}_3 \cup \left\{\overset{00}{\delta}\right\} \cup \{\alpha' \to \psi^7, \xi\}) = S_A]$,

where

$\overset{00}{\delta} = \psi^7 \to \overset{0}{\delta}, \xi = \psi^7 \to (\psi^1 \to \psi^{12})$,

$A^1 = Cn\left(R^P_{0+}, L^1_1 \cup {}_{\psi^7 O_0 u_{27} \sim\psi^1}L^r_2 \cup N\ddot{\overline{\Psi}}_3 \cup \{\xi\}\right)$.

Proof. Suppose to the contrary that

(1) $\langle R^P_{0+}, L^r_2 \cup A_r\rangle \notin Cns^T_A$

and

(2) $\neg(\exists\alpha' \in \overline{S}_A - A^1)\left(\forall \overset{00}{\delta} \in \overline{S}_A - A^1\right)$

$[Cn(R^P_{0+}, L^1_1 \cup {}_{\psi^7 O_0 u_{27} \sim\psi^1}L^r_2 \cup N\ddot{\overline{\Psi}}_3 \cup \left\{\overset{00}{\delta}\right\} \cup \{\alpha' \to \psi^7, \xi\}) = S_A]$,

where

(3) $\overset{00}{\delta} = \psi^7 \to \overset{0}{\delta}$,

(4) $\xi = \psi^7 \to (\psi^1 \to \psi^{12})$,

(5) $A^1 = Cn(R^P_{0+}, L^1_1 \cup {}_{\psi^7 O_0 u_{27} \sim\psi^1}L^r_2 \cup N\ddot{\overline{\Psi}}_3 \cup \{\xi\})$.

Hence,

(6) $(\forall\alpha \in \overline{S}_A - A^1)\left(\exists \overset{00}{\delta}{}' \in \overline{S}_A - A^1\right)$

$[Cn(R^P_{0+}, L^1_1 \cup {}_{\psi^7 O_0 u_{27} \sim\psi^1}L^r_2 \cup N\ddot{\overline{\Psi}}_3 \cup \left\{\overset{00}{\delta}{}'\right\} \cup \{\alpha \to \psi^7, \xi\}) = A^+ \neq S_A]$,

where

(7) $\overset{00}{\delta}{}' = \psi^7 \to \overset{0}{\delta}{}'$,

(8) $\xi = \psi^7 \to (\psi^1 \to \psi^{12})$,

(9) $A^1 = Cn(R^P_{0+}, L^1_1 \cup {}_{\psi^7 O_0 u_{27} \sim\psi^1}L^r_2 \cup N\ddot{\overline{\Psi}}_3 \cup \{\xi\})$.

Next, from $(I_{11})$, $(I_{14})$ and $(I_{16})$, it follows that

(10) $Cn(R^P_{0+}, L^1_T) \subseteq Cn(R^P_{0+}, L^1_1)$.

From (9) and $(I_{11})$, it follows that

(11) $\sim\psi^7 \in A^1 \Rightarrow \left(\forall\overset{0}{\delta} \in \overline{S}_A\right)\left[\psi^7 \to \overset{0}{\delta} \in A^1\right]$.

Hence, from (7), it follows that

(12) $\sim\psi^7 \in A^1 \Rightarrow \overset{00}{\delta}{}' \in A^1$.

Hence, from (6) – (9), it follows that

(13) $\sim\psi^7 \notin A^1$.

From (6) – (9), $(I_5)$, $(I_6) - (I_8)$, $(I_{14}) - (I_{17})$, it follows that

(14) $(\forall\alpha \in \overline{S}_A - A^1)\left(\exists \overset{00}{\delta}{}' \in \overline{S}_A - A^1\right)$

$[Cn(R^P_{0+}, L^1_1 \cup {}_{\psi^7 O_0 u_{27} \sim\psi^1}L^r_2 \cup N\ddot{\overline{\Psi}}_3 \cup \left\{\overset{00}{\delta}{}'\right\} \cup \{\alpha \to \psi^7, \xi\}) = A^+ \ \& \ A^+ \neq S_A \ \&$

$\gamma'_4, O_0 \to \gamma_0, \gamma'_2, O_0 \to \left(u_{27} \to \overset{1}{\beta}\right), O_0 \to \overset{1}{\beta}$,

$\gamma'_0 \to O_0, \gamma'_0 \to \overset{1}{\beta}, \overset{1}{\beta} \in A^+]$.

From (1), (6) – (10), (14), $(I_5) - (I_8)$, $(I_{11}) - (I_{17})$, by **Corollary 4.4.**, it follows that

(15) $(\forall\alpha \in \overline{S}_A - A^1)\left(\exists \overset{00}{\delta}{}' \in \overline{S}_A - A^1\right)$

$[Cn(R^P_{0+}, L^1_1 \cup {}_{\psi^7 O_0 u_{27} \sim\psi^1}L^r_2 \cup N\ddot{\overline{\Psi}}_3 \cup \left\{\overset{00}{\delta}{}'\right\} \cup \{\alpha \to \psi^7, \xi\}) = A^+ \ \& \ A^+ \neq S_A \ \& \ \psi^7 \to \overset{0}{\delta}{}'$,

$\psi^1 \to (\psi^7 \to \psi^{12}), \overset{1}{\beta}, \psi^{12} \wedge \psi^7 \wedge \psi^1 \to \overset{0}{\beta}$,

$\psi^{12} \to (\psi^7 \to \sim\psi^1), \psi^7 \to \sim\psi^1, \gamma'_4$,

$O_0 \to \gamma_0, \gamma'_2, O_0 \to \gamma'_0, \gamma'_0 \to O_0, \gamma'_0 \to u_{27}$,

$\psi^7 \to (\sim\psi^1 \to O_0), \psi^7 \to O_0, \psi^7 \to \gamma'_0$,

$\psi^7 \to u_{27} \in A^+]$.

From (15) and $(I_{12})$, it follows that

(16) $\left(\forall\alpha \in \overline{S}_A - A^1\right)\left(\exists\, \overset{00}{\delta}{}' \in \overline{S}_A - A^1\right)$

$\left[{}_{\psi^7 O_0 \sim \psi^1}L_2^r \subseteq A^+\right]$.

From (15) and (16), it follows that

(17) $\left(\forall\alpha \in \overline{S}_A - A^1\right)\left(\exists\, \overset{00}{\delta}{}' \in \overline{S}_A - A^1\right)[O_6 \in A^+]$,

where

(18) $O_6 = O_0 \to (\psi^7 \to \psi^1)$.

Hence, from (15) and $(I_6)$, it follows that

(19) $\left(\forall\alpha \in \overline{S}_A - A^1\right)\left(\exists\, \overset{00}{\delta}{}' \in \overline{S}_A - A^1\right)$

$[Cn(R_{0+}^P, L_1^1 \cup {}_{\psi^7 O_0 u_{27} \sim \psi^1}L_2^r \cup N\ddot{\overline{\Psi}}_3 \cup$

$\left\{\overset{00}{\delta}{}'\right\} \cup \{\alpha \to \psi^7, \xi\}) = A^+$ & $A^+ \neq S_A$ &

$\gamma_0' \to (\psi^7 \to \psi^1), \psi^7 \to \psi^1,$

$\psi^1 \to \sim\psi^7, \sim\psi^7 \in A^+]$.

Hence, from (9), it follows that

(20) $\left(\forall\alpha \in \overline{S}_A - A^1\right)\left(\exists\, \overset{00}{\delta}{}' \in \overline{S}_A - A^1\right)$

$[Cn(R_{0+}^P, L_1^1 \cup {}_{\psi^7 O_0 u_{27} \sim \psi^1}L_2^r \cup N\ddot{\overline{\Psi}}_3 \cup \left\{\overset{00}{\delta}{}'\right\} \cup$

$\{\alpha \to \psi^7, \xi\}) = A^+$ & $A^+ \neq S_A$ &

$A^+ - A^1 = \emptyset$ & $A^1 \subseteq A^+]$.

Hence, from (19), it follows that

(21) $\sim\psi^7 \in A^1$,

what contradicts with (13). □

## 5. The Main Result

**Theorem 5.1.** $\langle R_{0+}^P, L_2^r \cup A_r\rangle \in Cns_A^T$.

Proof. Let

1) $\langle R_{0+}^P, L_2^r \cup A_r\rangle \notin Cns_A^T$.

Hence, by **Lemma 4.3.** and **Lemma 4.4.**, we obtain that

2) $\left(\exists\alpha' \in \overline{S}_A - A^1\right)[Cn(R_{0+}^P, L_1^1 \cup$

${}_{\psi^7 O_0 u_{27} \sim \psi^1}L_2^r \cup N\ddot{\overline{\Psi}}_3 \cup \{\alpha' \to \psi^7\} \cup \{\xi\}) = S_A]$,

where

3) $A^1 = Cn(R_{0+}^P, L_1^1 \cup {}_{\psi^7 O_0 u_{27} \sim \psi^1}L_2^r \cup N\ddot{\overline{\Psi}}_3 \cup \{\xi\})$.

From 2) and 3), by **Theorem 3.5.**, it follows that

4) $\left(\exists\alpha' \in \overline{S}_A - A^1\right)[\sim(\alpha' \to \psi^7) \in A^1]$,

where

5) $A^1 = Cn(R_{0+}^P, L_1^1 \cup {}_{\psi^7 O_0 u_{27} \sim \psi^1}L_2^r \cup N\ddot{\overline{\Psi}}_3 \cup \{\xi\})$.

From 4) and 5), it follows that

6) $\left(\exists\alpha' \in \overline{S}_A\right)[\alpha' \notin A^1$ & $\alpha' \in A^1]$,

where

7) $A^1 = Cn(R_{0+}^P, L_1^1 \cup {}_{\psi^7 O_0 u_{27} \sim \psi^1}L_2^r \cup N\ddot{\overline{\Psi}}_3 \cup \{\xi\})$.

Contradiction □

**Theorem 5.2.** $\langle R_{0+}^P, L_2^{r'} \cup A_r'\rangle \in Cns_A^T$.

Proof. The proof of this **Theorem**, is analogical to the proof of **Theorem 5.1.** □